\documentclass[reqno,10pt]{amsart}    

\usepackage{amsmath}
\usepackage{amsfonts, amssymb}
\usepackage{mathrsfs,  bbm, cancel}
\usepackage[all]{xy}
\newtheorem{prop}{Proposition}[section]
\newtheorem{thm}[prop]{Theorem}
\newtheorem{lem}[prop]{Lemma}

\theoremstyle{definition}
\newtheorem{defn}[prop]{Definition}

\theoremstyle{remark}
\newtheorem{rem}[prop]{Remark}
\newtheorem{cor}[prop]{Corollary}

\newcommand{\z}{\mathbbm Z}
\renewcommand{\r}{\mathbbm R}

\newcommand{\X}{\mathfrak{X}}

\begin{document}
\title{Polynomial poly-vector fields}
\author{Frank Klinker}
\date{February 19, 2008}
\address{Faculty of Mathematics, Dortmund University of Technology, 44221 Dortmund, Germany} 
\address{frank.klinker@math.tu-dortmund.de}
\begin{abstract} In this text we give a decomposition result on polynomial poly-vector fields generalizing a result on the decomposition of homogeneous Poisson structures. 
We discuss consequences of this decomposition result in particular for low dimensions and low degrees. We provide the tools to calculate simple cubic Poisson structures in dimension three and quadratic Poisson structures in dimension four. 
Our decomposition result has a nice effect on the relation between Poisson structures and Jacobi structures.
\end{abstract}
\maketitle

\section{Introduction}

We motivate our discussion on  poly-vector fields by Poisson structures. For this we recall:
A Poisson structure
on the n-dimensional manifold $M$ is a skew symmetric bi-linear map $\{.,.\}:C^\infty(M)\times C^\infty(M)\to C^\infty(M)$
which obeys
\begin{itemize}
\item[P1:] The map $\{f,.\}:C^\infty(M)\to C^\infty(M)$ is   a derivation for all $f\in C^\infty(M)$.
\item[P2:] The Jacobi identity  $\{f,\{g,h\}\}+\{h,\{f,g\}\}+\{g,\{h,f\}\}=0$ holds for all $f,g,h\in C^\infty(M)$.
\end{itemize}
Due to P1 a Poisson structure is equivalent to a skew symmetric tensor field $\Pi\in\Gamma\Lambda^2 TM$. The Poisson bracket and the field are related by
$\{f,g\}=\Pi(df,dg)$. This may be generalized to brackets of degree higher than two, i.e.\ $\{\cdot,\ldots,\cdot\}:\times^{2l}C^\infty(M)\to C^\infty(M)$. If we impose conditions similar to  P1 and P2 we arrive at the so called generalized Poisson structures due to \cite{AzPerePerez2} and \cite{AzPerePerez}. 

To translate condition P2 to the tensor $\Pi$, we have to introduce  poly-vector fields.  
We denote the vector fields on $M$ by $\X^1(M)=\X(M)=\Gamma TM$
  and the $k$-vector fields by $\X^k(M):=\Lambda^k\X(M)$. A natural grading is given by $\X^k_{\text{nat}}(M):=\X^k(M)[-k]$.
The Lie bracket of vector fields extends to the Schouten bracket 
$[\;,\;]:\X^k(M)\times\X^l(M) \to \X^{k+l-1}(M)$ with
$[X,U]=L_XU$ for $X\in\X(M)$.
This turns the shifted algebra $\mathfrak{T}_{\text{poly}}(M) =\big({\bigoplus}_{k\geq 0} \X^k(M)[-k]\big)[1]$ into a
Gerstenhaber algebra. 
This means that the bracket obeys the graded Jacobi identity
\begin{equation}
(-)^{\tilde U\tilde W}[U,[V,W]]+(-)^{\tilde V\tilde W}[W,[U,V]]+(-)^{\tilde U\tilde V}[V,[W,U]]=0
\end{equation}
so that $\mathfrak{T}_{\text{poly}}(M)$ has a graded Lie algebra structure.
Furthermore, the bracket is  compatible with the degree one product, $\wedge :\X^k(M)\times\X^l(M)\to\X^{k+l}(M)$, on $\mathfrak{T}_{\text{poly}}(M) $   in the sense that 
\begin{equation}
[U,V\wedge W]=[U,V]\wedge W+(-)^{\tilde U(\tilde V+1)}V\wedge[U,W]\,
\end{equation}
The symmetry of the wedge product is given by $U\wedge V=(-)^{(\tilde U+1)(\tilde V+1)}V\wedge U$. 
To distinguish the shifted grading of $\mathfrak{T}_{\text{poly}}(M)$ from the natural grading  we add a tilde in the first case.
Turning back to the Poisson structures expressed by the associated bi-vector $\Pi$, condition P2 is connected  to the Schouten bracket. More precisely, a bi-vector field yields a Poisson structure if and only if  
\begin{itemize}
\item[P2':] $[\Pi,\Pi]=0$.
\end{itemize} 
In local coordinates, $\Pi=\Pi^{ij}\partial_i\wedge\partial_j$,  the latter is written as $\Pi^{[ij}(\partial_j\Pi^{k]l})=0$.

In the following we deal with polynomial poly-vector fields, i.e.\ poly-vector fields which coefficients are polynomials. 
For a  bi-vectors this means  we have a decomposition 
$\Pi=\sum \Pi^{ij}\partial_i\wedge\partial_j$ with polynomial coefficients $\Pi^{ij}$.
For two poly-vectors $\sigma_1$ and $\sigma_2$ which are homogeneous of degree $k$ and
$\ell$ the bracket $[\sigma_1,\sigma_2]$ is homogeneous of degree $k+\ell-1$. 
This yields that for a polynomial bi-vector field to obey [P2],  it is necessary that its highest homogeneous part
is a Poisson structure itself. For this reason, the restriction to homogeneous polynomial poly-vector fields
is natural. 

In the case of Poisson structures, the following holds. 
For degree zero all constant bi-vectors provide Poisson
structures. For degree one we write $\Pi^{ij}=f^{ij}_kx^k$ and condition [P2]
is given by $f^{[ij}_lf^{k]l}_n=0$, i.e.\ $f^{ij}_k$ are the structure
constants of a Lie algebra. For a classification of linear Poisson structures in four dimensions see \cite{Yunhe}.
For the discussion of quadratic Poisson structures we consider the injective  map $\imath_0: \mathfrak{gl}_n\to \X(M)$ with $\imath_0(A)=A^i{}_jx^j\partial_i \in\X(M)$ which gives the linear vector fields on $M$. 
The Lie bracket  on $\mathfrak{g}:=\mathfrak{gl}_n$ naturally extends to a Gerstenhaber bracket on $\widehat{\mathfrak{g}}=\big(\bigoplus_{k\geq0}\Lambda^k\mathfrak{g}[-k]\big)[1]$ and  $\imath_0$ to a morphism $\imath: \widehat{g}\to \mathfrak{T}_{\text{poly}}$ of Gerstenhaber algebras. A classical R-matrix is an element $r\in\Lambda^2\mathfrak{gl}_n$ which obey the Yang-Baxter equation $[r,r]=0$. A classical R-matrices yields a Poisson structure by $r=r^{ijkl}E_{ij}\wedge E_{kl}\mapsto \imath_1(r)=r^{ijkl}x^ix^k\partial_j\wedge\partial_l$. 
In \cite{BhaskRama} this connection between quadratic Poisson structures and classical R-matrices is recognized and a 1-1 correspondence is shown to be true at least in dimension $n=2$.  
In higher dimensions this correspondence is not true, because the graded components $\imath_k$  have non trivial kernels, in particular for $k=2,3$. For $n=3$ the authors in \cite{ManchonMasmoudiRoux} give an explicit example $r\in \widehat{\mathfrak{g}}^{(1)}= \Lambda^2\mathfrak{gl}_n$ such that $[r,r]\neq0$ and
$\imath_1(r)\neq0$ but $\imath_2[r,r]=[\imath_1(r),\imath_1(r)]=0$. Up to now the non linear Poisson structures resist descriptions as nice as in the linear case. Therefore, the only chance to describe them is to find an explicit way to construct them. 
This has been done for quadratic Poisson structures in \cite{DufourHaraki} and \cite{XuLiu} and used in \cite{ManchonMasmoudiRoux}.

In this text we give a decomposition formula for polynomial poly-vector fields into irreducible components, see Theorem \ref{maintheorem}. 
For any poly vector fields the Schouten bracket can be decomposed with respect to this decomposition, see Theorem \ref{bracketAB} and Corollary \ref{bracketAA}, and for even poly-vector fields the generalized Poisson condition P2' translates into compatibility conditions on the irreducible components, see Theorem \ref{poissonconditions} and Proposition \ref{ranktwo}. 
The decomposition result yields a connection between  polynomial Poisson structures and polynomial Jacobi structures which is collected in Theorem \ref{poissonjacobi}. The case where the polynomial degree and the vector field degree coincide turns out to be an exceptional combination. 
In the last part of this text we provide some tools to calculate polynomial Poisson structures for low dimensions and low degrees and give  a lot of examples.

\section{The trace operator}

We consider the volume form $\Psi$ on the oriented manifold $M$. $\Psi$ induces an isomorphism of poly-vector fields and differential forms. A $k$-vector field $U\in \X^k(M)$ is mapped to an ($n-k$)-form $\Psi(U)\in\Omega^{n-k}(M)$ by $ \Psi(U)(W) :=\Psi(U\wedge W)$.
In local coordinates with $\Psi=dx^1\wedge\cdots\wedge dx^n=\frac{1}{n!}\epsilon_{i_1\cdots i_n}dx^{i_1}\wedge\cdots\wedge dx^{i_n}$ and $U=\frac{1}{k!}U^{j_1\ldots j_k}\partial_{j_1}\wedge\cdots\wedge\partial_{j_k} $ this reads
\begin{equation}\label{localPsi}
\Psi(U)=\frac{1}{k!(n-k)!}\epsilon_{j_1\ldots j_k i_1\ldots i_{n-k}}U^{j_1\ldots j_k}dx^{i_1}\wedge\cdots\wedge dx^{i_{n-k}}\,.
\end{equation}
We use this isomorphism to define a derivation  $D:\mathfrak{T}_{\text{poly}}(M)\to\mathfrak{T}_{\text{poly}}(M)$ of degree $-1$ by
\begin{equation}\label{trace}
 D:  \X^k(M)  \to\X^{k-1}(M),\quad
D  =\Psi^{-1}\circ d\circ\Psi\,.
\end{equation}
This is in fact a differential, i.e.\ $D^2=0$. We note that  if we restrict $D$ to vector fields we recover the divergence with respect to the volume form $\Psi$. We will recall this soon. In local coordinates the action of $D$ is given by 
\begin{equation}\label{localD}
D(U)= \frac{1}{(k-1)!}\partial_m U^{i_1\cdots i_{k-1}m}\partial_{i_1}\wedge\cdots \wedge \partial_{i_{k-1}}\,,
\end{equation}
where we used $\epsilon^{i_1\ldots i_{\ell}a_1\ldots a_{n-\ell}}\epsilon_{j_1\ldots j_\ell a_1\ldots a_{n-\ell}}=\ell!(n-\ell)! \delta^{i_1\ldots i_{\ell}}_{j_1\ldots j_\ell} $.

With respect to the non-shifted grading of poly-vector fields, $D$ is compatible with the wedge product and the Schouten bracket in  the following sense. 
\begin{prop}
For all poly-vectors $U$ and $U'$ we have 
\begin{equation}
D(U\wedge U')= (-)^{U'}DU\wedge U'  + U\wedge DU'+(-)^{U'}  [U,U']\,. \label{compatible}
\end{equation}
\end{prop}

\begin{proof}
Because both sides are additive we only have to check the identity for homogeneous $U$ and $U'$. Consider local coordinates with $\Psi=dx^1\wedge\ldots\wedge dx^n$ and let $U=f\partial_{i_1}\wedge\cdots\wedge\partial_{i_{k}}$ and $U'=g\partial_{i_1}\wedge\cdots\wedge\partial_{i_{k'}}$. We evaluate the four terms of (\ref{compatible}) separately and use the short notation $\partial_{i_1}\wedge\cdots\wedge\partial_{i_k}= \partial_{i_1\cdots i_k}$.
The left-hand side  is given by 
\begin{align*}
D(U\wedge U')
=\ & D(fg\,\partial_{i_1\cdots\wedge i_k j_1\cdots j_{k'}}) \\
=\ & (-)^{k+k'-1}(k+k')\partial_{[i_1}(fg)\partial_{i_2\cdots i_k j_1\cdots j_{k'}]}\\
=\ &  (-)^{k+k'-1} k \partial_{[i_1}(fg)\partial_{i_2 \cdots i_k]j_1 \cdots j_{k'}}
	+(-)^{k'-1}k'\partial_{[j_1}(fg)\partial_{|i_1 \cdots i_k| j_2\cdots j_{k'}]}\\
=\ &  (-)^{k+k'-1} k f\partial_{[i_1}g\,\partial_{i_2\cdots i_k] j_1\cdots j_{k'}}
 	+  (-)^{k+k'-1} k g\partial_{[i_1}f\,\partial_{i_2 \cdots i_k] j_1\cdots j_{k'}}\\
& 	+(-)^{k'-1}k'f\partial_{[j_1}g\,\partial_{|i_1 \cdots i_k| j_2 \cdots  j_{k'}]}
	+(-)^{k'-1}k'g\partial_{[j_1}f\,\partial_{|i_1 \cdots i_k| j_2 \cdots\ j_{k'}]}\,.
\end{align*}
The summands of the right-hand side are
\begin{align*}
U\wedge DU' 
=\ &  (-)^{k'-1}k'f \,\partial_{i_1 \cdots i_k}\wedge \partial_{[j_{1}}g\,\partial_{j_2 \cdots j_{k'}]}\\
=\ & (-)^{k'-1}k'f\partial_{[j_{1}}g\, \partial_{|i_1 \cdots i_k | j_2 \cdots j_{k'}]}\,, \\
\end{align*}
and
\begin{align*}
DU\wedge U' 
=\ & (-)^{k-1}k \partial_{[i_1}f\, \partial_{i_2\cdots i_k]}\wedge 
		g \partial_{j_1 \cdots j_{k'}}\\
=\ & (-)^{k-1}k g \partial_{[i_1}f\, \partial_{i_2\cdots i_k] j_1 \cdots j_{k'}}\,,
\end{align*}
as well as 
\begin{align*}
 \big[U,\ U'\big] 
 =\ &      	f\big[\partial_{i_1 \cdots\ i_k},g\big]\wedge \partial_{j_1 \cdots j_{k'}} \\
& 	-(-)^{(k-1)(k'-1)} g \big[\partial_{j_1 \cdots j_{k'}},f\big] \wedge\partial_{i_1 \cdots i_k} \\
=\ & 	f \big(\sum_{\alpha=1}^{k}(-)^{k-\alpha} \partial_{i_\alpha}g\,
	 	\partial_{i_1\cdots\widehat{i_\alpha} \cdots i_k}\big) \wedge \partial_{j_1 \cdots j_{k'}} \\
& 	-(-)^{(k-1)(k'-1)} g\big(  \sum_{\alpha=1}^{k'} (-)^{k'-\alpha}\partial_{j_{\alpha}}f\,
		\partial_{j_1 \cdots \widehat{ j_{\alpha} } \cdots  j_{k'}}\big)\wedge\partial_{i_1\cdots i_k} \\
=\ & 	(-)^{k-1} f \big(\sum_{\alpha=1}^{k}(-)^{\alpha-1} \partial_{i_\alpha}g\, 
	 	\partial_{i_1\cdots\widehat{i_\alpha}\cdots i_k}\big) \wedge \partial_{j_1 \cdots j_{k'}} \\
& 	+(-)^{kk'+k-1} g\big(  \sum_{\alpha=1}^{k'} (-)^{\alpha-1}\partial_{j_{\alpha}}f\,
		\partial_{j_1 \cdots \widehat{ j_{\alpha}} \cdots j_{k'}}\big)  \wedge\partial_{i_1 \cdots i_k} \\
=\ & 	(-)^{k-1} k f \partial_{[i_1}g\, \partial_{i_2 \cdots i_k] j_1 \cdots j_{k'}} 
	+(-)^{kk'+k-1}k' g \partial_{[j_1}f\,\partial_{j_2\cdots j_{k'}] i_1 \cdots i_k}\\
=\ & 	(-)^{k-1} k f \partial_{[i_1}g\, \partial_{i_2 \cdots i_k] j_1 \cdots j_{k'}} 
	+(-)^{-1} k' g \partial_{[j_1}f\,\partial_{|i_1 \cdots i_k| j_2 \cdots j_{k'}]}\,.
\end{align*}
When we compare these terms we see that (\ref{compatible}) holds.
\end{proof}

In the next lemma  we collect some relations between $[\cdot,\cdot]$ and $D$. 
\begin{lem}\label{UDU}
Let $U$ and $V$ be any poly vector fields. Then 
\begin{align}
 [DU,U]&=-[U,DU],\quad DU\wedge U= U\wedge DU \label{eins1}\\
\intertext{and}
D[U,V]&=[U,DV]-(-)^V[DU,V]\,. \label{zwei} 
\end{align}
In particular,
for $U$ even  $D[U,U]=-2[DU,U]$, and so $[DU,U]=0$ if $[U,U]=0$.
\end{lem}

\begin{proof}
To prove (\ref{eins1}) we recall the commutation properties of the wedge product and the Schouten bracket and get $[DU,U]=-(-)^{(\tilde U-1)\tilde U}[U,DU]=-[U,DU]$ as well as $DU\wedge U=(-)^{(U-1)U}U\wedge DU=U\wedge DU$. 
\begin{align*}
D^2(U\wedge V)  =\ 
    & D\big( (-)^V DU\wedge V+ U\wedge DV +(-)^V [U,V] \big) \\
=\ & (-)^V DU\wedge DV +(-)^V(-)^V[DU,V]+(-)^{V-1}DU\wedge DV \\ 
    &+(-)^{V-1}[U,DV]+(-)^VD[U,V]\\
=\ & (-)^V\big( D[U,V]-[U,DV]+(-)^V[DU,V]\big)\,.
\end{align*}
$D^2=0$ yields
\[
D[U,V]=[U,DV]-(-)^V[DU,V] =[U,DV]-(-)^{(U-1)(V-1)}[V,DU]
\]
which proves (\ref{zwei}).  
With $U=V$ even  
this is $D [U,U] = -2[DU,U]$.
\end{proof}

We mentioned before  that if we restrict $D$ to vector fields we recover the divergence with respect to the volume form $\Psi$. In this case  (\ref{compatible})  reads as  
$D(X\wedge Y)= ({\rm div}Y)X- ({\rm div}X)Y  - [X,Y]$.
A connection to the notion of divergence induced by a connection $\nabla$ on the manifold $M$ is given by the next Proposition.
\begin{prop}{\cite{KobNom1}}
Let $\Psi$ be a volume form on the manifold $M$. Consider a connection $\nabla$ on $M$ with vanishing  torsion $T=0$ which is compatible with the volume form, i.e.\ $\nabla\Psi=0$. Then the two notions of divergence of vector fields given by
\begin{equation}
{\rm div}_\Psi X:=\Psi^{-1}\circ d\circ\Psi(X)\quad\text{and}\quad {\rm div}_\nabla X:=\text{tr}(\nabla X)
\end{equation}
coincide.
\end{prop}

The existence of a differential $D$ does not depend on whether the manifold is orientable or not, although in our case it will be. This is due to 
\begin{prop}{\cite{Koszul1}}\ \ 
Let $M$ be a manifold with torsion free connection $\nabla$. Then there exists a differential  $D_\nabla:\X^k(M)\to\X^{k-1}(M)$ with (\ref{compatible}) and which coincides  for vector fields $X$ with the divergence associated to $\nabla$.
\end{prop}

We will next concentrate on polynomial poly-vector fields $\mathfrak{P}\subset\mathfrak{T}_{\text{poly}}$. These are poly-vector fields whose coefficients in a local coordinate system are polynomials in the coordinates. 
The space $\mathfrak{P}$ is bi-graded, because every subspace $\mathfrak{P}^{(k)}\subset \mathfrak{T}_{\text{poly}}^{(k-1)}$ admits a further $\z$-grading due to the degree of its coefficient polynomial. 
We denote these spaces by $\mathfrak{P}^{(\ell)}=\bigoplus_k \mathfrak{P}^{(k,\ell)}$.   
An element $A\in\mathfrak{P}^{(k,\ell)}$ can be written as $A=\frac{1}{k!\ell!}A_{i_1\ldots i_k}{}^{j_1\ldots j_\ell}x^{i_1}\cdots x^{i_k}\partial_{j_1}\wedge\cdots\wedge \partial_{j_\ell}$ with  $A_{\bullet}{}^{\bullet}$ totally symmetric in its lower indices and totally skew symmetric in its upper indices.
\begin{defn}\label{homogen}
A polynomial $\ell$-vector field $A$ is called $k$-homogeneous if $A\in\mathfrak{P}^{(k,\ell)}$.
\end{defn}
$A_\bullet{}^\bullet$  transforms like a tensor with the given index picture, if we restrict to  linear coordinate transformations. In this case we identify
\[
\mathfrak{P}^{(k,\ell)}=S^k(\Omega(M))\otimes\Lambda^\ell(\X(M))\,.
\]
The restriction of the differential on $\mathfrak{T}_{\text{poly}}$ to the polynomial poly-vector fields is part of the next Proposition.
\begin{prop}
The differential $D$ on $\mathfrak{T}_{\text{poly}}$, cf.\  (\ref{trace}), restricted to the polynomial poly-vector fields $\mathfrak{P}$ is the trace operator. In particular, with respect to the bi-grading of $\mathfrak{P}$ the  differential $D$ is of degree $(-1,-1)$. I.e.
\begin{equation}\begin{split}
&D: S^k(\Omega(M))\otimes\Lambda^\ell(\X(M))\to S^{k-1}(\Omega(M))\otimes\Lambda^{\ell-1}(\X(M)) \\
& D(A)_{i_1\ldots i_{k-1}}{}^{j_1\ldots j_{\ell-1}}= A_{ i_1\ldots i_{k-1}m}{}^{ j_1\ldots j_{\ell-1}m}
\end{split}\end{equation}
\end{prop}

\begin{cor}
For  linear vector fields $A$, i.e.\ for elements in $\Omega(M)\otimes\X(M)=\text{End}(\X(M))$, the differential $D$ is the usual trace for endomorphisms, i.e.\ $D(A)= \text{tr} A$.
\end{cor}
\begin{proof}
We write the poly-vector field in the form 
\[
A=\frac{1}{k!\ell!}A_{i_1\ldots i_k}{}^{j_1\ldots j_\ell}\,x^{i_1}\cdots x^{i_k}
\partial_{j_1}\wedge\cdots\wedge \partial_{j_\ell}\,.
\]
We recall (\ref{localD}) which yields
\begin{align*}
D(A)
=\ & \tfrac{1}{(\ell-1)!k!} \partial_m \big(A_{i_1\ldots i_k}{}^{j_1\ldots j_{\ell-1}m}\,x^{i_1}\cdots x^{i_k}\big) 
\partial_{j_1}\wedge\cdots\wedge \partial_{j_{\ell-1}}\\
=\ & \tfrac{1}{(\ell-1)!k!} A_{i_1\ldots i_k}{}^{j_1\ldots j_{\ell-1}m}\partial_m \big(x^{i_1}\cdots x^{i_k}\big) 
\partial_{j_1}\wedge\cdots\wedge \partial_{j_{\ell-1}}\\
=\ & \tfrac{1}{(\ell-1)!k!} A_{i_1\ldots i_k}{}^{j_1\ldots j_{\ell-1}m}
\sum_{\alpha=1}^k \delta_m^{i_\alpha} x^{i_1}\cdots \widehat{x^{i_\alpha}}\cdots x^{i_k}
\partial_{j_1}\wedge\cdots\wedge \partial_{j_{\ell-1}}\\
=\ & \tfrac{1}{(\ell-1)!(k-1)!} A_{i_1\ldots i_{k-1} m}{}^{j_1\ldots j_{\ell-1}m}
\,x^{i_1}\cdots x^{i_{k-1}}
\partial_{j_1}\wedge\cdots\wedge \partial_{j_{\ell-1}}\,.
\end{align*}
\end{proof}

\section{The decomposition of poly-vector fields}

Let $V$ be an $n$ dimensional real vector space. We decompose  $\mathfrak{P}^{(k,\ell)}= S^k(V^*)\otimes \Lambda^\ell V$  into its irreducible representation spaces with respect to $\mathfrak{sl}_n\r$. 
\begin{defn}
We denote the representation space to the weight $\lambda$ by $V_\lambda$. Furthermore, we denote by $\lambda_m$ the weight $(\underbrace{0,\ldots,0,1}_m)$ and for $\lambda=k\lambda_{n-1}+\lambda_{\ell}$ we use the short notation $V_{k\ell}$.
\end{defn}
We recall the isomorphism of representations
\[
V^*\simeq \Lambda^{n-1}V=V_{\lambda_{n-1}}\,.
\]
The space $S^k(V^*)\simeq S^k(\Lambda^{n-1}V)$ is itself irreducible and its weight is given by $\lambda=k\lambda_{n-1}=(0,\ldots,0,k)$. To get the decomposition of $\mathfrak{P}^{(k,\ell)}$ we calculate 
\begin{align}
\mathfrak{P}^{(k,\ell)}
&=S^k(V^*)\otimes\Lambda^\ell V\simeq V_{k\lambda_{n-1}}\otimes V_{\lambda_\ell} \nonumber\\
&= V_{k\lambda_{n-1}+\lambda_\ell}\oplus V_{\lambda_n+(k-1)\lambda_{n-1}+\lambda_{\ell-1}}\nonumber \\
&\simeq  V_{k\lambda_{n-1}+\lambda_\ell}\oplus V_{(k-1)\lambda_{n-1}+\lambda_{\ell-1}}\nonumber \\
&= V_{k,\ell}\oplus V_{k-1,\ell-1}\label{decomposition}\,.
\end{align}
The dimension of the representation space $V_{k,\ell}$ is
\[
\text{dim}V_{k,\ell}=\frac{(n+k)!}{(n+k-\ell)k!\ell!(n-\ell-1)!}\, ,
\]
in particular, $\text{dim}V_{k,k}=\frac{1}{(k!)^2}\prod_{j=1}^{k}(n^2-j^2) $.  

We note that $V_{1,1}=\mathfrak{sl}_n$, $\mathfrak{P}^{(k,0)}=S^k(V^*)=V_{k\lambda_{n-1}}=V_{k,0} $, $\mathfrak{P}^{(0,\ell)}=\Lambda^\ell(V)=V_{\lambda_{\ell}}=V_{0,\ell} $ and  $\mathfrak{P}^{k,n}\simeq S^k(V^*)=V_{k-1,n-1}$.

The trace operator $D: \mathfrak{P}^{(k,\ell)}\to \mathfrak{P}^{(k-1,\ell-1)}$ gives rise to an exact sequence. For $m>0$ it is 
\begin{align}\label{sequence}
S^{n+m}(V^*) &\overset{D}{\longrightarrow}
S^{n+m}(V^*)\oplus V_{n+m-2,n-2}\overset{D}{\longrightarrow} 
\cdots  
\nonumber\\
\cdots &				\overset{D}{\longrightarrow}
V_{m+2,2}\oplus V_{m+1,1}	\overset{D}{\longrightarrow}
V_{m+1,1}\oplus S^m(V^*)		\overset{D}{\longrightarrow}
S^m(V^*)\,, \\
\intertext{or for $0<m<n$ }
S^{n-m}(V^*)&			\overset{D}{\longrightarrow}
S^{n-m}(V^*)\oplus V_{n-m-2,n-2}\overset{D}{\longrightarrow} 
\cdots  
\nonumber\\
\cdots &				\overset{D}{\longrightarrow}
V_{2,m+2}\oplus V_{1,m+1}	\overset{D}{\longrightarrow}
V_{1,m+1}\oplus \Lambda^m V	\overset{D}{\longrightarrow}
\Lambda^m V \,,\\
\intertext{or for $m=0$}
S^{n}(V^*)&			\overset{D}{\longrightarrow}
S^{n}(V^*)\oplus V_{n-2,n-2} 	\overset{D}{\longrightarrow} 
\cdots  
\nonumber\\
\cdots&				\overset{D}{\longrightarrow}
V_{3,3}\oplus V_{2,2}		\overset{D}{\longrightarrow}
V_{2,2}\oplus \mathfrak{sl}_n	\overset{D}{\longrightarrow}
\mathfrak{sl}_n \oplus \r		\overset{D}{\longrightarrow}
\r\,.
\end{align}

The kernel and the image of the $\mathfrak{sl}_n$-invariant trace are representation spaces, too. Due to the irreducibility of the decomposition (\ref{decomposition}) this yields
\begin{equation}
\text{ker}(D|_{\mathfrak{P}^{(k,\ell)}})=V_{k,\ell},\qquad 
\text{im}(D|_{\mathfrak{P}^{(k,\ell)}})=V_{k-1,\ell-1}\,.
\end{equation}

We  recall the wedge product 
\begin{align*}
\wedge : \mathfrak{P}^{(k,\ell)}\otimes\mathfrak{P}^{(k'\ell')} & \ \longrightarrow \mathfrak{P}^{(k+k',\ell+\ell')}
\end{align*} 
to describe the inclusion
\begin{equation}
S^{k-1}(V^*)\otimes\Lambda^{\ell-1}V \supset V_{k-1,\ell-1} 
\hookrightarrow
V_{k-1,\ell-1}\subset S^k(V^*)\otimes\Lambda^\ell V\,. 
\end{equation}
Consider the element $e^{(k,\ell)}\in\mathfrak{P}^{(1,1)}$ defined by 
\begin{equation}
e^{(k,\ell)}=\frac{1}{n+k-\ell}\,e_0, \quad\text{with } e_0=x^m\partial_m \,.
\end{equation}  
This vector field has certain properties collected in the next lemma.
\begin{lem}
 For $A\in\mathfrak{P}^{(k',\ell')}$ we have
\begin{gather}
D(e^{(k,\ell)})=\frac{n}{n+k-\ell}\,, \\
\big[e^{(k,\ell)},A\big]\ =\ \frac{k'-\ell'}{n+k-\ell} A\,. \label{liee}
\end{gather}
If furthermore $DA=0$, we have
\begin{equation}
D(A\wedge e^{(k,\ell)})=\frac{n+k'-\ell'}{n+k-\ell}A\,.
\end{equation}
\end{lem}
\begin{proof}
The trace of $e_0$  is $n$, so the first part is obvious. The second part follows  from 
$[e_0,fY] =[e_0 ,f] Y+f[e_0,Y]$ and 
\begin{align*}
\big[ e_0, x^{i_1}\cdots x^{i_{k'}}\big] 
&=  x^m \partial_m (x^{i_1}\cdots x^{i_{k'}})\\
&= x^m \sum_{\alpha=1}^{k'} \delta_m^{i_\alpha} x^{i_1}\cdots \widehat{x^{i_{\alpha}}}\cdots x^{i_{k'}}\\
&=k' x^{i_1}\cdots x^{i_{k'}}\,,
\end{align*}
as well as
\begin{align*}
\big[ e_0, \partial_{j_1}\wedge\cdots\wedge\partial_{j_{\ell'}}\big] 
&= \partial_m \wedge \big[ x^m, \partial_{j_1}\wedge\cdots\wedge\partial_{j_{\ell'}}\big] \\
&= \partial_m \wedge \sum_{\alpha=1}^{\ell'} (-)^\alpha \delta_m^{j_\alpha} 
       \partial_{j_1}\wedge\cdots\wedge\widehat{\partial_{j_{\alpha}}}\wedge\cdots\wedge\partial_{j_{\ell'}}\\
&= -\ell' \partial_{j_1}\wedge\cdots\wedge\partial_{j_{\ell'}}\,.
\end{align*}
We  plug $A$ and $e^{(k,\ell)}$ into (\ref{compatible}) and  with $DA=0$ as well as (\ref{liee}) we get
\begin{align*}
D(A\wedge e^{(k,\ell)}) 
&=  A\wedge De^{(k,\ell)}-\big[ A,e^{(k,\ell)}\big]\\
&= \frac{n}{n+k-\ell}A+\big[ e^{(k,\ell)},A\big] \\
&= \ \frac{n+k'-\ell'}{n+k-\ell}A \,.
\end{align*}
\end{proof}

\begin{cor}
The preimage of an element $A\in V_{k-1,\ell-1}\subset S^{k-1}(V^*)\otimes\Lambda^{\ell-1}V$ in 
$V_{k-1,\ell-1}\subset S^k(V^*)\otimes\Lambda^\ell V$ with respect to $D$ is given by
\[
A\wedge e^{(k,\ell)}\,.
\]
The vector $e^{(k,\ell)}$ does only depend on the difference $(k-\ell)$ and so is fixed for the sequence (\ref{sequence}).
\end{cor}
This discussion proves the following theorem.
\begin{thm}\label{maintheorem}
Each polynomial poly-vector field $A\in\mathfrak{P}^{(k,\ell)}$ admits a unique decomposition $A=A_0+A_1$ with $A_0\in V_{k,\ell}\subset\mathfrak{P}^{(k,\ell)} $ and $A_1\in V_{k-1,\ell-1}\subset\mathfrak{P}^{(k,\ell)}$. Explicitly the fields are connected by $DA=DA_1$, $DA_0=0$, $A_1=DA\wedge e^{(k,\ell)}$, i.e.
\begin{equation}\label{eins}
A=A_0+DA\wedge e^{(k,\ell)}\, .
\end{equation}
\end{thm}

In the next step we write the bracket of  two polynomial poly vector field in terms of the poly vector fields we get by applying the decomposition with respect to (\ref{eins}).
\begin{thm}\label{bracketAB}
Let $A\in\mathfrak{P}^{(k,\ell)}$ and $B\in\mathfrak{P}^{(k',\ell')}$ polynomial poly vector fields 
with decompositions $A=A_0+DA\wedge e^{(k,\ell)}$ and  $B=B_0+DB\wedge e^{(k',\ell')}$ cf.\ (\ref{eins}). Then   

\begin{equation}\begin{split}
\big[A,B\big]_0 =\ &  [A_0,B_0\big] \\
& + \frac{\Delta'}{n+\Delta}DA\wedge B_0
	+(-)^{\ell'}\frac{\Delta}{n+\Delta'}A_0\wedge DB \\
& +\frac{\Delta'}{n+\Delta}\big[A_0,DB\big]\wedge e^{(k'',\ell'')}  
	-(-)^{\ell'}\frac{\Delta}{n+\Delta'} \big[DA,B_0\big]\wedge e^{(k'',\ell'')} 
\end{split}\end{equation}

and
\begin{equation}\begin{split}
D\big[A,B\big] =\ &  \big[A_0,DB\big]-(-)^{\ell'}\big[DA,B_0\big] \\
&	- \frac{(n+\Delta'')(\Delta-\Delta')}{(n+\Delta)(n+\Delta')}\Big(DA\wedge DB 
		+(-)^{\ell'}\big[DA,DB\big]\wedge e^{(k'',\ell'')}\Big)\,,
\end{split}\end{equation}
with $\Delta=k-\ell$, $k''=k+k'$, $\ell''=\ell+\ell'$, $\Delta'=k'-\ell'$ and $\Delta''=\Delta+\Delta'=k''-\ell''$.
\end{thm}

\begin{proof}
Let $A$ and $B$ as above as stated. We write $e=e^{(k,\ell)}$, $ e'=e^{(k',\ell')}$, and $ e''= e^{(k'',\ell'')}$. 
From the symmetry properties of the wedge product and the bracket and from the compatibility (\ref{compatible}) we get
\begin{align}
&[A\wedge V,B\wedge Y] \nonumber\\ 
=\ &	(-)^{(\tilde A+1)(\tilde V+\tilde B+1)}[V,B]\wedge A\wedge Y 
	+(-)^{\tilde V\tilde Y+\tilde B(\tilde Y+1)+1}[A,Y]\wedge V\wedge B\nonumber\\
&	+(-)^{\tilde B(\tilde V+1)}[A,B]\wedge V\wedge Y
	+(-)^{\tilde A(\tilde Y+\tilde V)+\tilde B(\tilde Y+1)+\tilde V+\tilde A}[V,Y]\wedge A\wedge B\,.
	\label{wedgebracket}
\end{align}
For $V=e$ and $Y=e'$ this is 
\begin{equation}
[A\wedge e,B\wedge e']=
A\wedge [e,B]\wedge e' - [e',A]\wedge B\wedge e\,. \label{wedgebracket2} 
\end{equation}
We make use of the decompositions of $A$ and $B$  and formulas  (\ref{wedgebracket2}) as well as (\ref{liee}) and get
\begin{align*}
\big[A,B\big]=\ & \big[A_0+DA\wedge e, B_0+DB\wedge e' \big]\\
=\ &   \big[A_0, B_0 \big]
 + \big[A_0,DB\wedge e' \big]
 + \big[DA\wedge e, B_0\big] 
 + \big[DA\wedge e, DB\wedge e' \big] \displaybreak[3]\\
=\ & \big[A_0, B_0 \big]
 + \big[A_0,DB\big] \wedge e'  +(-)^{(\ell-1)(\ell'-1)}DB\wedge [A_0,e']  \\
 & + DA\wedge \big[e, B_0\big] +(-)^{\ell'-1}\big[DA,B_0\big]\wedge e \\
 & +DA\wedge [e,DB]\wedge e' - [e',DA]\wedge DB\wedge e \displaybreak[3]\\
=\ & \big[A_0, B_0 \big] 
 + (-)^{\ell'} \frac{\Delta}{n+\Delta'}A_0\wedge DB +\frac{\Delta'}{n+\Delta}DA\wedge B_0\\
& +\frac{1}{n+\Delta'} \big[A_0,DB\big] \wedge e_0  
	-(-)^{\ell'}\frac{1}{n+\Delta}\big[DA,B_0\big]\wedge e_0 \\
& +  \frac{\Delta'-\Delta}{(n+\Delta)(n+\Delta')}  DA\wedge DB\wedge e_0\,.
\end{align*}
To get the trace of this expression we use lemma \ref{UDU}.
\begin{align*}
& D\big[A,B\big] = \big[A,DB\big]-(-)^{\ell'}\big[DA,B\big] \\
 =\ & \big[A_0,DB\big]-(-)^{\ell'}\big[DA,B_0\big] 
         +\big[DA\wedge e,DB\big]-(-)^{\ell'}\big[DA,DB\wedge e'\big] \displaybreak[3]\\
 =\ & \big[A_0,DB\big]-(-)^{\ell'}\big[DA,B_0\big]
         +DA\wedge\big[ e,DB\big]+(-)^{\ell'}\big[DA,DB\big]\wedge e \\
    &  -(-)^{\ell'}\big[DA,DB\big]\wedge e'-(-)^{\ell'+\ell(\ell'-1)}DB\wedge \big[DA,e'\big]\displaybreak[3]\\
=\ &  \big[A_0,DB\big]-(-)^{\ell'}\big[DA,B_0\big] 
     +\Big(\frac{\Delta'}{n+\Delta}-\frac{\Delta}{n+\Delta'} \Big)DA\wedge DB \\
    & + (-)^{\ell'}\Big(  \frac{1}{n+\Delta}-\frac{1}{n+\Delta'} \Big) \big[DA,DB\big]\wedge e_0 \,.
\end{align*}
Using $\Delta(n+\Delta)-\Delta'(n+\Delta')=(n+\Delta'')(\Delta-\Delta')$ proves the second part of Theorem \ref{bracketAB}. 
The first part is obtained by rearranging the summands in  $\big[A,B]_0= \big[A,B] - D[A,B]\wedge e''$.
\end{proof}

If we restrict to the case $A=B$ we  get the following consequences.

\begin{cor}\label{bracketAA}
Let $A\in\mathfrak{P}^{(k,\ell)}$ with $\ell$ even. The bracket of $A$ with itself in terms of its decomposition with respect to (\ref{eins}) is given by
\begin{align}
\big[ A,A\big]_0 	
       &= \big[ A_0,A_0 \big] 
	+ \frac{2(k-\ell)}{n+k-\ell}\big(DA \wedge A_0 -\big[DA,A_0\big]\wedge e^{(2k,2\ell)} \big) \,,
	\label{bracketAA1} \\
D\big[A,A\big] 	
       &=-2[DA,A_0]\,.\label{bracketAA2}
\end{align}
In $[DA,A_0]$ and $DA \wedge A_0$ we may replace $A_0$ by $A$.
\end{cor}
\begin{cor}\label{coreasy}
All polynomial poly-vector fields of type $A\wedge e^{(k,\ell)}$ have vanishing bracket $[A\wedge e^{(k,\ell)},A\wedge e^{(k,\ell)}]$.
\end{cor}

From Corollary \ref{bracketAA} we immediately  conclude the following result on even homogeneous polynomial poly vector fields, e.g.\ (generalized) Poisson structures. 

\begin{thm}\label{poissonconditions} 
Let $A\in\mathfrak{P}^{(k,\ell)}$ and $\ell$ even with $A=A_0+DA\wedge e^{(k,\ell)}$ cf.\ Theorem \ref{maintheorem}.  
\begin{enumerate}
\item  Suppose $k\neq\ell$. Then $\big[A,A\big]=0$ if and only if  
\begin{equation}\label{pcond1}
\big[A_0,A_0\big]=\frac{2(\ell-k)}{n+k-\ell} DA\wedge A_0\,.
\end{equation} 
\item For $k=\ell$ we have $\big[A,A\big]=0$ if and only if  
\begin{equation}\label{pcond2}
\big[A_0,A_0\big]=\big[DA,A_0\big]=0\,. 
\end{equation} 
\end{enumerate}
\end{thm} 

\begin{rem}
We emphasize on the fact that in the first point, the vanishing of $[DA,A_0]$ is hidden in the stated equation (\ref{pcond1}). This is seen by taking the trace on both sides of (\ref{pcond1})  and then using (\ref{zwei}) as well as (\ref{compatible}).

The further specialization to $k=\ell=2$ or $\ell=2$ gives the decomposition result on quadratic Poisson structures as stated in \cite{XuLiu} or \cite{MalekShafei}.
\end{rem}

At the end of this section we will give a  remark on diffeomorphisms and  their effect  on the decomposition cf.\ Theorem \ref{maintheorem}. 
A diffeomorphism which respects the polynomial structure as well as the degree  has to be linear. 
For a linear diffeomorphism $L(x)^j=L_i{}^jx^i$ and a polynomial poly-vector field we write $L_*A$ for the induced action. E.g.\  for a linear vector field $A\in\mathfrak{gl}_n$ we have $L_*A=L^{-1}AL$, in particular $L_*e^{(k,\ell)}=e^{(k,\ell)}$.
The isomorphism $\Psi$ cf.\ (\ref{localPsi}) is compatible with this diffeomorphism in the sense that 
for poly-vector fields $A$ and $B$ we have
\begin{equation}\label{LundPsi}
L^*\circ \Psi\circ L_*A (B)=\Psi(L_*A,L_*B)=\det L\, \Psi(A,B)=\det L\,\Psi(A) (B)\,
\end{equation}
or $L^*\circ\Psi\circ L_*=\det L\, \Psi$. Furthermore, $d$ commutes with $L^*$ so that we get
\begin{equation}
DL_*A=\Psi^{-1}d\Psi L_*A=\det L\, L_*\Psi^{-1}d\Psi A=\det L\, L_*DA\,.
\end{equation}
This yields
\begin{prop}\label{compatibility}
Let $A,\hat A\in\mathfrak{P}^{(k,\ell)}$ be  polynomial poly-vector fields with trace decomposition $A=A_0+B\wedge e$ and $\hat A=\hat A_0+\hat B\wedge e$, respectively. Then $A$ and $\hat A$ are diffeomorphic if and only if there exist a linear diffeomorphism $L$ such that $A_0=L_*\hat A_0$ and $B=L_* \hat B$.
\end{prop}

\section{Applications to Poisson and Jacobi structures}
\subsection{Poisson structures}

Theorem \ref{poissonconditions} has a nice consequence in the case $k\geq 3$. Let $A\in\mathfrak{P}^{(k,2)}$ be a Poisson structure with non vanishing trace. Furthermore suppose that its trace free part $A_0$ is  a Poisson structure, too.  Due to the vanishing of the expression $DA\wedge A_0$, the trace free part has the form $A_0=DA\wedge C$ for a vector field $C$.  This motivates the next definition.
\begin{defn}\label{trivialPoisson}
A (generalized) Poisson structure is called simple if its trace-free part is itself a (generalized) Poisson structure. 
\end{defn}
In the case of Poisson structures we may summarize
\begin{prop}\label{ranktwo}
A simple homogeneous polynomial Poisson structure of degree $\geq 3$  has rank two. 
\end{prop}

{\bf Bi-vector fields in dimension two. }  
For $V=\r^2$ we have $V^*\simeq V$ as representation spaces and $V_{m-1,1}\simeq S^{m}(V)$. In particular, $[A,A]=0$ for all $A\in\mathfrak{P}^{(m,2)}(V)$. The sequence (\ref{sequence}) is written as
\[
S^m(V)\overset{D}{\longrightarrow} S^{m}(V)\oplus S^{m-2}(V)\longrightarrow S^{m}(V)
\]
and we have $\text{ker}D=\{0\}$.

Let us consider $m=2$. We  have $\mathfrak{sl}_2\simeq S^2(V)$ via 
\[
\begin{pmatrix}a &b\\ c & -a\end{pmatrix} \mapsto \begin{pmatrix} -b&a \\ a& c  \end{pmatrix}
\]
which is the $(\partial_1\wedge\partial_2)$-component of $\frac{1}{2}\begin{pmatrix}a &b\\ c & -a\end{pmatrix}\wedge e^{(1,1)}$.
This yields the result mentioned in the introduction, namely that every quadratic Poisson structure is obtained by an element in $\Lambda^2\mathfrak{gl}_2$. This, of course, may be generalized to arbitrary $m\geq2$. 
\begin{prop}
Every homogeneous polynomial Poisson structure of degree $m$ in dimension two  is obtained by 
\[
S^{m-1}(V^*)\otimes V\supset S^m(V) \ni A\mapsto A\wedge e^{(m,2)}\in S^m(V^*)\otimes\Lambda^2(V)\,.
\]
\end{prop}

{\bf Linear bi-vector fields. }
Consider linear bi-vectors and the splitting $\mathfrak{P}^{(1,2)}=V_{1,2}\oplus \Lambda^1V$.
\begin{itemize}
\item  
The bi-vector with trace $v\in\Lambda^1V\simeq V$ is $A=\frac{1}{n-1}v\wedge\mathbbm{1}$, i.e.\ $A^{ij}{}_k=\frac{2}{n-1}v^{[i}\delta^{j]}_k$. 
These coefficients give rise to a Lie algebra ($\hat=$ Poisson structure  in the linear  case) due to Corollary \ref{coreasy}. 

\item 
In dimension $n=3$ this may be expressed by the help of  the cross product. Let $E=(E_1,E_2,E_3)$ be a basis of the Lie algebra. The bracket can be expressed as follows
\begin{align*}
\big[E_i,E_j\big] =\ & v^{[i}\delta^{j]}_kE_k = \tfrac{1}{2}\big(v^{i}\delta^{j}_kE_k-v^{j}\delta^{i}_kE_k\big) \\ 
=\ & \tfrac{1}{2}\big(v^{i}E_j-v^{j}E_i\big)=\tfrac{1}{4}\epsilon^{ijk}(v\times E)_k
\end{align*}
E.g. $v=2e_1$ gives 
\[
\big[E_1,E_2\big]=E_2,\quad  \big[E_1,E_3\big]=E_3,\quad  \big[E_2,E_3\big]=0
\]
which is the Lie algebra $\mathfrak{g}_{(1,1)}:=\r\ltimes \r^2$ of rank 2.
Furthermore, consider the Lie algebra $\mathfrak{g}_{(\alpha,\beta)}$ given by
\[
\big[E_1,E_2\big]=\alpha E_2,\quad  \big[E_1,E_3\big]=\beta E_3,\quad  \big[E_2,E_3\big]=0
\]
This algebra has structure constants $A^{12}{}_2=\alpha, A^{13}{}_3=\beta$ and trace $DA=(\alpha+\beta)e_1$. 
The decomposition (\ref{eins}) together with the discussion above yields 
\begin{gather*}
(A_0)^{12}{}_2=\tfrac{\alpha-\beta}{2}, (A_0)^{13}{}_3=\tfrac{\beta-\alpha}{2} \\
(A_1)^{12}{}_2 = (A_1)^{13}{}_3=\tfrac{\alpha+\beta}{2}\,.
\end{gather*}
In particular $[A_0,A_0]= DA\wedge A_0=0$, i.e.\ we have a decomposition of $\mathfrak{g}_{(\alpha,\beta)}$ into a pair of Lie algebras  $\big(\mathfrak{g}_{(\frac{\alpha-\beta}{2},\frac{\beta-\alpha}{2})}, \mathfrak{g}_{(\frac{\alpha+\beta}{2},\frac{\alpha+\beta}{2})}\big)$.

\item
Consider a semi simple Lie algebra and its decomposition $\mathfrak{g}=\mathfrak{h}\oplus\bigoplus_\alpha\mathfrak{g}_\alpha$. The roots $\alpha\in\mathfrak{h}^*$ obey $[H,X]=\alpha(H)X$ for all $H$ in the Cartan algebra $\mathfrak{h}$ and $X$ in the one dimensional subspace $\mathfrak{g}_\alpha$. The $\alpha(H)$ are the only non-vanishing diagonal elements of the structure constants and with $\alpha(H)$ the element $-\alpha(H)$ is present, too. So the sum over these elements vanish and so does the trace of $\mathfrak{g}$ which is given by
\[
D\mathfrak{g}= \sum_{i} \big(\sum_\alpha \alpha(H_i)\big)H_i
\]
for a basis $\{H_i\}$ of $\mathfrak{h}$.
\item 
Let $\mathfrak{g}=\mathfrak{g}_0\oplus\mathfrak{g}_1$ be the vector space decomposition of the Lie algebra $\mathfrak{g}$ into its semi simple and solvable parts with basis $\{E_i\}$ and $\{F_j\}$, respectively. We have $[\mathfrak{g}_0,\mathfrak{g}_1]\subset\mathfrak{g}_1$ so that the trace of $\mathfrak{g}$ decomposes like
\begin{equation}
D\mathfrak{g}
 =D\mathfrak{g}_1 + \sum_i \text{tr}(\text{ad}E_i|_{\mathfrak{g}_1}) E_i\,.
\end{equation}
A similar but not so strong decomposition of the trace may be obtained by the decomposition $\mathfrak{g}=\mathfrak{g}^{(1)}\oplus\mathfrak{g}/\mathfrak{g}^{(1)}$ with $\mathfrak{g}^{(1)}=[\mathfrak{g},\mathfrak{g}]$ and $[\mathfrak{g}/\mathfrak{g}^{(1)},\mathfrak{g}^{(1)}]\subset\mathfrak{g}^{(1)}$. It is not so strong in the sense that $\mathfrak{g}_1\subset \mathfrak{g}/\mathfrak{g}^{(1)}$ and $\mathfrak{g}_0\supset\mathfrak{g}^{(1)}$.
\end{itemize}

\subsection{Jacobi structures}

A (generalized) Jacobi structure on the manifold $M$ is by definition a pair $(\Lambda,E)\in\X^{2\ell}(M)\times\X^{2\ell-1}(M)$ which obeys
\begin{equation}\label{jacobi}
\big[\Lambda,E\big]=0
\quad\text{and}\quad 
\big[\Lambda,\Lambda\big]=2E\wedge\Lambda\,.
\end{equation}
This is due to \cite{Perezbueno} but with the field $E$ rescaled.
\begin{defn}
A Jacobi structure $(\Lambda,E)$ is called $k$-homogeneous  if $\Lambda\in\mathfrak{P}^{(k,2\ell)}$ and $E\in\mathfrak{P}^{(k-1,2\ell-1)}$ are homogeneous of degree $k$ and $k-1$, respectively.  

A $k$-homogeneous (generalized) Poisson structure $\Pi$ and a $k$-homogeneous (generalized) Jacobi structure $(\Lambda, E)$ are called associated, if 
\begin{equation}
\Pi_0=\Lambda_0\quad\text{and}\quad E_0= \frac{k-2\ell}{n+k-2\ell}(D\Lambda-D\Pi)\,.
\end{equation}
\end{defn}
\begin{thm}\label{poissonjacobi}
\begin{enumerate}
\item
For every homogeneous Jacobi structure $(\Lambda,E)$ of degree $k\neq 2\ell$ a $k$-homo\-geneous Poisson structure is defined 
via its decomposition with respect to Theorem \ref{maintheorem}
\[
\Pi_0=\Lambda_0
\quad\text{and}\quad
D\Pi=D\Lambda+\frac{n+k-2\ell}{2\ell-k}E_0\,.
\]
\item
Let $\Pi=\Pi_0+D\Pi\wedge e^{(k,2\ell)}$ be a $k$-homogeneous Poisson structure, $k\neq2\ell$, which admits a splitting $D\Pi=F_0+\tilde E_0$ such that there exists $\xi\in\mathfrak{P}^{(k-2,2\ell-2)}$ with
\[
\xi\wedge \Pi_0 \ =\ [\Pi_0 ,F_0]+ F_0\wedge \tilde E_0 \,.
\]
This defines a Jacobi structure given by
\[
\Lambda= \Pi_0+F_0\wedge e^{(k,2\ell)}\quad\text{and}\quad E= \frac{2\ell-k}{n+k-2\ell}\big(\tilde E_0+\xi\wedge e^{(k,2\ell)}\big)\, .
\]
\item
In the case $k=2\ell$ a pair of associated Jacobi- and  Poisson structure obeys $E_0=0$ and $D\Lambda=D\Pi+\eta$ for some $\eta$ with $[\eta,\Lambda_0]=-DE\wedge \Lambda_0$.
\end{enumerate}
\end{thm}

\begin{rem}
The second item of the preceding theorem includes two special cases, namely two structures coming from $\xi=0$.
The first is the obvious one with $E=0$, i.e. $F_0=D\Lambda$, which yields the Jacobi structure
\[
\Lambda=\Pi,\qquad E=0\,.
\]
The second special case  is $F_0=0$ with Jacobi structure given by
\[
\Lambda=\Pi_0,\qquad E=\frac{2\ell-k}{n+k-2\ell}D\Pi\,.
\]
For $\ell=1$ the latter is the structure also given in \cite{Petalidou}.
\end{rem}

\begin{proof}

Let $\Lambda=\Lambda_0+D\Lambda\wedge e$ and $E=E_0+DE\wedge e$ form a $k$-homogeneous Jacobi structure with decompositions with respect to Theorem \ref{maintheorem}. We use the  notations $e=e^{(k,2\ell)}$, $e''=e^{(2k,4\ell)}$, and $\Delta=k-2\ell$.  
We reformulate the Jacobi conditions (\ref{jacobi}) in terms of the fields $\Lambda_0$, $E_0$, $D\Lambda$ and $DE$ by using Theorem \ref{bracketAB}. We plug in $A=E,B=\Lambda$, $\ell'=2\ell$, and $\Delta=\Delta'=k-2\ell$. This yields 
\begin{align*}
\big[E,\Lambda\big]_0  =\ &  \big[E_0,\Lambda_0\big] \\
&	+\frac{\Delta}{n+\Delta}\Big( E_0\wedge D\Lambda +DE\wedge\Lambda_0
	+\big( [E_0,D\Lambda]-[DE,\Lambda_0]\big)\wedge e''\Big) 
\end{align*}
and $D[E,\Lambda]=[E_0,D\Lambda]-[DE,\Lambda_0]=- D(E_0\wedge D\Lambda+ DE\wedge\Lambda_0)$.
So the first Jacobi condition is written as
\[
\big[E,\Lambda\big]=0  \Longleftrightarrow 
\left\{ \begin{matrix} \displaystyle
\big[E_0,\Lambda_0\big] =-\frac{\Delta}{n+\Delta}\big(E_0\wedge D\Lambda +DE\wedge\Lambda_0\big) 
& \quad\text{if }\Delta\neq0\,,
    \\[2ex]
\big[E_0,\Lambda_0\big]=0\ \text{ and }\ \big[E_0,D\Lambda\big]-\big[DE,\Lambda_0\big]=0 
& \quad\text{if }\Delta=0\,.
\end{matrix}\right.
\]
We use (\ref{bracketAA}) for $A=\Lambda$ as well as
$E\wedge \Lambda= E_0\wedge \Lambda_0+\big(DE\wedge \Lambda_0+E_0\wedge D\Lambda\big)\wedge e$
and 
$D(E\wedge \Lambda)=DE\wedge \Lambda_0 +E_0\wedge D\Lambda  +[E,\Lambda]$. This yields for the second Jacobi condition
\begin{align*}
&\big[\Lambda,\Lambda\big] =\ 2E\wedge\Lambda\\ 
&\Longleftrightarrow 
\left\{\begin{array}{l}
	\displaystyle  
	\big[\Lambda_0,\Lambda_0\big]
	+\frac{2\Delta}{n+\Delta}D\Lambda\wedge \Lambda_0 
	-\frac{2\Delta}{n+\Delta}\big[D\Lambda,\Lambda_0\big]\wedge e'' 
\\[1.5ex]
	\displaystyle \quad 
	= 2 E_0\wedge \Lambda_0 
	+\Big(\frac{2\Delta}{n+\Delta}\big( DE\wedge \Lambda_0+E_0\wedge D\Lambda\big)
	-2\big[E,\Lambda\big]\Big)\wedge e'' 
\\[1.5ex] \qquad\quad \text{and} \\[1.5ex]
	\big[D\Lambda,\Lambda_0\big] =  -DE\wedge \Lambda_0 -E_0\wedge D\Lambda  -[E,\Lambda] 
\end{array}\right.
\\
&\Longleftrightarrow 
\left\{\begin{array}{l}
	\displaystyle   
	\big[\Lambda_0,\Lambda_0\big] 
	= 2 ( E_0- \frac{\Delta}{n+\Delta}D\Lambda) \wedge \Lambda_0 
	-2\frac{n+2\Delta}{n+\Delta} \big[E,\Lambda\big]\wedge e''
\\[1.5ex]\qquad\quad \text{and} \\[1.5ex]
	\big[D\Lambda,\Lambda_0\big]= -DE\wedge \Lambda_0 -E_0\wedge D\Lambda  -\big[E,\Lambda\big]\,.
\end{array}\right.
\end{align*}
So the two Jacobi conditions together translate as follows
\begin{align*}
& \big[E,\Lambda\big]=0 \quad\text{and}\quad \big[\Lambda,\Lambda\big ]=2E\wedge\Lambda 
\\
&\Longleftrightarrow 
\left\{
\begin{array}{rlc}
  \displaystyle 
	\big[D\Lambda, \Lambda_0\big] 
	= &  -DE\wedge \Lambda_0 -E_0\wedge D\Lambda 
		&\qquad(*)
\\[1.5ex]\text{and}&&\\[1.5ex]
 \displaystyle
	\big[\Lambda_0,\Lambda_0\big] 
	= & \displaystyle   2 ( E_0- \frac{\Delta}{n+\Delta}D\Lambda) \wedge \Lambda_0 
		&\qquad(**)
\end{array} 
\right.
\end{align*}
Suppose $\Delta\neq 0$, i.e. $k\neq2\ell$. Then for $\Pi_0:=\Lambda_0$ and $D\Pi= D\Lambda-\frac{n+\Delta}{\Delta}E_0$ we have
\[
\big[\Pi_0,\Pi_0\big]=-\frac{2\Delta}{n+\Delta}D\Pi\wedge \Pi_0\, ,
\]
i.e.\ (\ref{pcond1}). This proves the first item. 
The second item is almost the same up to the fact that we have to make sure the compatibility $(*)$ of $\Lambda$ and $E_0$ with the trace $DE=\xi$. 
In the third item the vanishing of $E_0$ follows from the definition. Therefore, $D\Pi$ and $D\Lambda$ are not related via $E$. Then the  Jacobi structure $(\Lambda=\Lambda_0+D\Lambda\wedge e, E=\xi\wedge e)$ and the Poisson structure $\Pi=\Lambda_0+D\Pi\wedge e$ are associated only if the difference $\eta= D\Lambda-D\Pi $ obeys $[\eta,\Lambda_0]=-\xi\wedge \Lambda_0$, i.e.\ $D\Pi$ must coincide with some part of $D\Lambda$ which leaves $\Lambda_0$ invariant.
\end{proof}

\section{Tools for explicit calculations} 

In this section we provide the tools to give  the lists of simple cubic Poisson structures in dimension three and  quadratic Poisson structures in dimension four. The latter supplements the list of quadratic Poisson structures in dimension three  which may be found e.g.\  in \cite{XuLiu} or \cite{DufourHaraki}.

\subsection{Cubic Poisson structures in dimension three}\label{cubic3}

Let $\Pi$ be a cubic Poisson structure. Theorem \ref{maintheorem} provides a decomposition of the bi-vector field $\Pi$ into a pair $(\Pi_0,A_0)$. 
Here $\Pi_0$ denotes the trace free part of $\Pi$ and depends on a biquadratic polynomial $f$, via the isomorphism (\ref{localPsi}). 
The vector field $A_0=A-DA\wedge e^{(3,2)}$ is the trace free part of a quadratic vector field $A$ and encodes the trace of $\Pi$. 

The Poisson condition $[\Pi,\Pi]=0$  may be translated to the function $f$ and the vector field $A_0$. 
We do this for the equivalent conditions  $[\Pi_0,\Pi_0]=A_0\wedge \Pi_0$ and $[A_0,\Pi_0]=0$, compare Theorem \ref{poissonconditions}. 
We recall  the calculations in  (\ref{liederivation}) and that  $[\Pi_0,\Pi_0]=D(\Pi_0\wedge\Pi_0)$ vanish  in  dimension three, i.e.\ $\Pi$ is simple. 
So we are left with the conditions 
\begin{equation}\label{poissonfA}
A_0(f)=0\quad\text{and}\quad A_0\wedge \Pi_0=0\,.
\end{equation}
For $\Pi_0$ being a two-form, the last equality forces it to be of the form $\Pi_0=A_0\wedge C$ for a  vector field $C$  (compare proposition \ref{ranktwo}) and the Poisson structure is of course of rank $\leq 2$. If we translate the Poisson condition with respect to these data, we get
\begin{equation}\label{poissonCA}
L_CA_0=0
\end{equation}
where we used $[A_0,A_0]=0$. For a trace free vector field $C$ the operators $D$ and $L_C$ commute so that the condition may be read $L_C A=0$.

To gain a list of all simple Poisson structures with trace free vector field $C\in \mathfrak{P}^{(1,1)}$ we restrict to its  Jordan form due to proposition \ref{compatibility}. 
We write the components of the vector field $A=A^k\partial_k$ with respect to the basis of quadratic polynomials, i.e.\ $A^k=A_{11}^kx^2+A_{22}^ky^2+A_{33}^kz^2+A_{12}^kxy+A_{13}^kxz+A_{23}^kyz$ or 
\begin{equation*}
A^k= \begin{pmatrix}
A_{11}^k,A_{22}^k,A_{33}^k,A_{12}^k,A_{13}^k, A_{23}^k
\end{pmatrix}^T
\end{equation*}
The bracket of $A$ with the vector field $C=C_i^jx^i\partial_j$ is given by
\begin{equation*}
[C,A]^k=C_i^jx^i\frac{\partial A^k}{\partial x^j}-C^k_jA^j\,.
\end{equation*}
Because of the restriction to the Jordan form of $C$, we will only need the following vectors for the calculations
\begin{multline*}
\big(
x\partial_xA^k , y\partial_yA^k, z\partial_zA^k, x\partial_yA^k,y\partial_zA^k,y\partial_xA^k
\big) = \\
\begin{pmatrix}
\begin{matrix}
2A_{11}^k\\ 0\\0\\A^k_{12}\\A_{13}^k\\0
\end{matrix}\ \ 
\begin{matrix}
0\\2A_{22}^k\\ 0\\A^k_{12}\\0\\A_{23}^k
\end{matrix}\ \ 
\begin{matrix}
0\\0\\2A_{33}^k\\0\\A^k_{13}\\A_{23}^k
\end{matrix}\ \ 
\begin{matrix}
A_{12}^k\\0\\0\\2A^k_{22}\\A_{23}^k\\0
\end{matrix}\ \ 
\begin{matrix}
0\\A_{23}^k\\0\\A^k_{13}\\0\\2A_{33}^k
\end{matrix}\ \ 
\begin{matrix}
0\\A_{12}^k\\0\\2A^k_{11}\\0\\A_{13}^k
\end{matrix}
\end{pmatrix}
\end{multline*}

The next list gives  all pairs  of trace less $(C,A_0) \in V_{1,1}\times V_{2,1}\subset \mathfrak{P}^{(1,1)}\times \mathfrak{P}^{(2,1)}$ which yield a simple Poisson structure 
\[
\Pi=A_0\wedge\big(C+e^{(3,2)}\big)
\] 
in the described way. This list is complete up to a permuting the variables.

\begin{itemize}
\item[{[A]}] 
	$C=\begin{pmatrix}a_1&&\\&a_2&&\\&&a_3\end{pmatrix}$ with $a_1+a_2+a_3=0$.
\begin{itemize}
\item[{[A.1]}] 
	$a_1a_2a_3\neq 0$.
\begin{itemize}
\item[{[A.1.1]}] 
	$a_i\neq 2a_j$ for all $1\leq i,j\leq 3$. 

	In this case $A=0$ is the only solution of $[C,A]=0$.
\item[{[A.1.2]}]  
	$a_2=2a_1$, i.e.\ $C=a\cdot \text{diag}(1,2,-3)$.

	We have one  non-zero solution of $[C,A]=0$ and this solution is trace free. It is given by 
	\begin{equation*}
		A_0=\alpha x^2\partial_y\,.
	\end{equation*}
\end{itemize}
\item[{[A.2]}] 
	$a_1=0,a_2\neq0$.

	The space of solutions $A$ is six-dimensional and we find a trace free basis 
	\begin{equation*}
		{\renewcommand{\arraystretch}{1.5} 
		\begin{array}{r l} 
		A_0\in\text{span}\big\{ 
		\!\!\!\!\! 	& y^2\partial_x, \   z^2\partial_x,\   yz\partial_x,\  x^2\partial_x-xy\partial_y-xz\partial_z, \\ 
		\!\!\!\!\!	& x^2\partial_x-3xy\partial_y+xz\partial_z,\  x^2\partial_x+xy\partial_y-3xz\partial_z\,  
		\big\}\,. 
		\end{array}}
	\end{equation*}
\item[{[A.3]}] 
	$a_1=a_2=0$.

	In this the case the Poisson structure is of the form $\Pi=A_0\wedge e^{(3,2)}$ with an arbitrary trace free quadratic vector field $A_0$.
\end{itemize}
\item[{[B]}] 
	$C=\begin{pmatrix}a&1&\\&a&\\&&-2a \end{pmatrix}$.
\begin{itemize}
\item[{[B.1]}]
	$a\neq0$. 

	There is no non-zero solution of $[C,A]=0$.
\item[{[B.2]}] 
	$a=0$. 

	The space of solutions is eight dimensional and the projection onto the trace free part is six dimensional:
	\begin{equation*}
		{\renewcommand{\arraystretch}{1.5}
		\begin{array}{rl}
		A_0\in\text{span}\big\{ 
		\!\!\!\!\!	&z^2\partial_z-(xz\partial_x+yz\partial_y), \ 3xz\partial_z-(x^2\partial_x+xy\partial_y), \\
		\!\!\!\!\!	&  x^2\partial_y,\  z^2\partial_y,\  xz\partial_y, \  x^2\partial_z\,  \big\}\,.
		\end{array}}
	\end{equation*}
\end{itemize}
\item[{[C]}] 
	$C=\begin{pmatrix}0&1&\\&0&1\\&&0 \end{pmatrix}$.

	The dimension of the space of solutions of $[C,A]=0$ is four. The projection on the trace free part is three dimensional and given by
	\begin{equation*}
		A_0\in\text{span}\big\{\,
			x^2\partial_z,\  xy\partial_z-x^2\partial_y,\ 
			xy\partial_y+ x^2\partial_x+2y^2\partial_z-3xz\partial_z\, 
		\big\}
	\end{equation*}
\item[{[D]}] 
	$C=\begin{pmatrix}a&b&\\-b&a&\\&&-2a \end{pmatrix}$ with $b\neq 0$. 
\begin{itemize}
\item[{[D.1]}] 
	$a\neq 0$. 

	As before we have no solution in this generic case.
\item[{[D.2]}] 
	$a=0$. 

	The solutions of $[C,A]=0$ span a three dimensional subspace of $\mathfrak{P}^{(2,1)}$ which has a trace free basis. 
	So the Poisson structures are built up by
	\begin{equation*}
		A_0\in\text{span} \big\{\, (y^2+x^2)\partial_z, \  z(y\partial_x-x\partial_y),\  z(x\partial_x+y\partial_y)\,  \big\}\,.
	\end{equation*}
\end{itemize}
\end{itemize}

\subsection{Quadratic Poisson structures in dimension four}\label{QuadraticPoisson4}

In this particular dimension we handle with quadratic bi-vectors and  their trace given by linear vector fields. 
A Poisson structure $\Pi$ is described by a pair $(\Pi_\theta,A)$  where $\Pi_\theta$ is a trace free Poisson structure and $A$ is linear vector field which is compatible with $\Pi_\theta$ in the way that $L_A\Pi_\theta=[A,\Pi_\theta]=0$ holds, cf.\ Theorem \ref{poissonconditions}.
The linear vector field encodes the trace of $\Pi$ and so is trace free, i.e.\ $A\in\mathfrak{sl}_4$. 
Furthermore, because $\Pi_\theta$ is trace free, there exists a cubic tri-vector field $L$ such that $\Pi_\theta=DL$. 
If we consider the isomorphism (\ref{localPsi})  we may write $L=\Psi^{-1}\theta$ with a cubic one-form $\theta$.
The trace free bi-vector is then given by $\Pi_\theta=\Psi^{-1}d\theta$, cf.\ (\ref{trace}). This explains the notation for the trace free part of $\Pi$.
The compatibility condition may be reformulated by
\begin{equation}\label{liederivation}\begin{split}
[A,\Pi_\theta] & = D(A\wedge\Pi_\theta) =\Psi^{-1}\circ d\circ \Psi(A\wedge\Pi_\theta) 
   = \Psi^{-1}\circ d\circ \imath_A\circ \imath_{\Pi_\theta}\Psi \\ 
& =\Psi^{-1}\circ d\circ \imath_A d\theta 
 =\Psi^{-1}(-\imath_A\circ d^2\theta+ L_A  d\theta) 
 =\Psi^{-1}\circ L_A d\theta \\
&=\Psi^{-1}\circ d L_A\theta = \Pi_{L_A\theta}
\end{split}\end{equation}
and so reads as 
\begin{enumerate}
\item[(i)]  $L_A\theta=0$\,.
\end{enumerate} 

The Poisson condition on $\Pi_\theta$ translates to $\theta$ in the following way. The commutator of $\Pi_\theta$ with itself is a tri-vector field and via $\Psi$ a one-form. We look at its components 
\begin{equation*}
(\Psi\circ [\Pi_\theta,\Pi_\theta])_m
=  (\Psi\circ D(\Pi_\theta\wedge\Pi_\theta))_m= \partial_m  \Psi(\Pi_\theta\wedge\Pi_\theta) \,.
\end{equation*} 
The function  $\Psi(\Pi_\theta\wedge\Pi_\theta)$ which is the only coefficient in $d\theta\wedge d\theta$ is a homogeneous  bi-quadratic polynomial. Its derivative vanishes if and only of it vanishes itself. So the Poisson condition on $\Pi_\theta$ can be written as
\begin{enumerate}
\item[(ii)] $d\theta\wedge d\theta=0$.
\end{enumerate}
To characterize all quadratic Poisson structures in the four-dimensional case we have to look for pairs of one-forms $\theta$ with cubic coefficients and trace free matrices $A$ such that (i) and (ii) are satisfied. We write $\Pi_{\theta,A}$ for the resulting Poisson structure.
We call a one-form and a matrix  which fulfills condition (i) compatible.  This is a linear condition on the coefficients of the 1-form. Because of remark \ref{propo1} below, which is an easy consequence of proposition \ref{compatibility} and (\ref{LundPsi}) as well as (i) and (ii),  we may restrict ourself to matrices $A$ which are in Jordan form. Condition (ii) is non-linear and yields algebraic relations of degree two for the coefficients.
\begin{rem}\label{propo1}
$\Pi_{\theta,A}$ and $\Pi_{\eta,B}$ are Poisson isomorphic if and only if there is a linear isomorphism $L$ such that 
\begin{equation*}
\eta=\det L\, L^*\theta,\quad\text{and}\ B=LAL^{-1}\,.
\end{equation*}
\end{rem}

For the explicit calculations in (i) we expand  the cubic one-form in the form $\theta=\theta_kdx^k$ and 
\begin{equation}
\theta_k = \sum_{0\leq m\leq n\leq o\leq 3} \theta_{k;mno}x^mx^nx^o \label{order}
\end{equation}
with $(x^0,x^1,x^2,x^3):=(t,x,y,z)$.
For $A\in\mathfrak{sl}_4$ we have
\begin{equation}\label{lieder}
L_A\theta_k=A^i{}_k\theta_i+A^i{}_jx^j\frac{\partial \theta_k}{\partial x^i}
\end{equation}
and $L_A\theta=0$ is a system on the 80 coefficients $\theta_{k;mno}$.
We write the coefficients of the respective 1-forms as vectors in the basis of the cubic polynomials as given in (\ref{order}).
In view of (\ref{lieder}) and the restriction to the Jordan form of the matrix $A$ we only need the following vectors:
\begin{gather*}
\big (
\theta_k,\ 
t\partial_t \theta_k,\  x\partial_x \theta_k,\  y\partial_y\theta_k,\  z\partial_z \theta_k ,\
t\partial_x\theta_k,\  x\partial_y\theta_k,\  y\partial_z\theta_k,\  x\partial_t\theta_k,\  z\partial_y\theta_k   
\big)=\\[1ex]
\begin{pmatrix}
\begin{matrix}
\theta_{k;012} \\\theta_{k;013}\\  \theta_{k;023}\\   \theta_{k;123}\\  \theta_{k;001}\\  \theta_{k;002}\\ \theta_{k;003}\\
\theta_{k;110}\\  \theta_{k;112}\\  \theta_{k;113}\\  \theta_{k;220}\\  \theta_{k;221}\\ \theta_{k;223}\\
\theta_{k;330}\\ \theta_{k;331}\\  \theta_{k;332}\\  \theta_{k;000}\\  \theta_{k;111}\\  \theta_{k;222}\\ \theta_{k;333}
\end{matrix}
\ 
\begin{matrix}
\theta_{k;012}\\  \theta_{k;013}\\ \theta_{k;023}\\ 0\\ 2\theta_{k;001}\\ 2\theta_{k;002}\\  2\theta_{k;003}\\  
\theta_{k;110}\\  0\\ 0\\ \theta_{k;220}\\  0\\  0\\ \theta_{k;330}\\  0\\  0\\ 3\theta_{k;000}\\  0\\  0\\  0
\end{matrix}
\  
\begin{matrix}
\theta_{k;012} \\  \theta_{k;013} \\  0 \\  \theta_{k;123} \\ \theta_{k;001} \\ 0 \\  0 \\  2\theta_{k;110} \\  2\theta_{k;112} \\ 
2\theta_{k;113} \\  0 \\  \theta_{k;221} \\  0 \\ 0 \\  \theta_{k;331} \\  0 \\ 0 \\  3\theta_{k;111} \\  0 \\ 0
\end{matrix}
\  
\begin{matrix}
\theta_{k;012} \\ 0\\ \theta_{k;023} \\ \theta_{k;123} \\0 \\ \theta_{k;002} \\ 0 \\ 0\\ \theta_{k;112} \\ 0\\ 2\theta_{k;220} \\ 
2\theta_{k;221} \\ 2\theta_{k;223} \\0 \\ 0 \\ \theta_{k;332} \\0\\ 0 \\ 3\theta_{k;222} \\ 0
\end{matrix}
\  
\begin{matrix}
0 \\ \theta_{k;013} \\ \theta_{k;023} \\ \theta_{k;123} \\0 \\ 0 \\ \theta_{k;003} \\ 0\\ 0\\ \theta_{k;113} \\0 \\ 0\\ \theta_{k;223} \\
2\theta_{k;330} \\ 2\theta_{k;331} \\ 2\theta_{k;332} \\0\\ 0 \\ 0\\ 3\theta_{k;333}
\end{matrix}
\ 
\begin{matrix}
2\theta_{k;112} \\ 2\theta_{k;113} \\ \theta_{k;123} \\ 0 \\ 2\theta_{k;110} \\ \theta_{k;012} \\\theta_{k;013} \\3\theta_{k;111}\\
0\\0\\\theta_{k;221}\\0\\0\\\theta_{k;331}\\0\\0\\\theta_{k;001}\\0\\0\\0
\end{matrix}
\  
\begin{matrix}
2\theta_{k;220}\\\theta_{k;023}\\0\\2\theta_{k;223}\\\theta_{k;002}\\0\\0\\\theta_{k;012}\\2\theta_{k;221}\\\theta_{k;123}\\
0\\3\theta_{k;222}\\0\\0\\\theta_{k;332}\\0\\0\\\theta_{k;112}\\0\\0
\end{matrix}
\  
\begin{matrix}
\theta_{k;013}\\0\\2\theta_{k;330}\\2\theta_{k;331}\\0\\\theta_{k;003}\\0\\0\\\theta_{k;113}\\0\\\theta_{k;023}\\\theta_{k;123}\\
2\theta_{k;332}\\0\\0\\3\theta_{k;333}\\0\\0\\\theta_{k;223}\\0
\end{matrix}
\  
\begin{matrix}
2\theta_{k;002} \\ 2\theta_{k;003} \\ 0\\ \theta_{k;023} \\3\theta_{k;000} \\ 0 \\ 0 \\ 2\theta_{k;001} \\ \theta_{k;012} \\ 
\theta_{k;013} \\0 \\ \theta_{k;220} \\ 0 \\0 \\ \theta_{k;330} \\ 0 \\0 \\ \theta_{k;110} \\ 0 \\ 0
\end{matrix}
\  
\begin{matrix}
0 \\ \theta_{k;012} \\ 2\theta_{k;220} \\ 2\theta_{k;221} \\0 \\ 0\\ \theta_{k;002} \\ 0 \\ 0 \\ \theta_{k;112} \\0 \\ 0 \\ 
3\theta_{k;222} \\\theta_{k;023} \\ \theta_{k;123} \\ 2\theta_{k;223} \\0\\ 0\\ 0 \\ \theta_{k;332}
\end{matrix}
\end{pmatrix}
\end{gather*}
A careful examination yields 8 cases with 43 subcases in total. We will restrict here to three examples. The whole list may be found in \cite{KlinkerWeb}.
\begin{itemize}
\item $A=\text{diag}(a_1,a_2,a_3,a_4)$, $|a_i|$ distinct, and no relation of the form $a_i=-3a_j$ holds for any $i,j$.

The only 1-form $\theta$ for which (i) holds is 
\[ 
\theta=\alpha_0xyz\,dt+\alpha_1tyz\,dx+\alpha_2txz\,dy+\alpha_3txy\,dz\,.
\] 
In particular, (ii) is also satisfied. the associated Poisson structure $\Pi_\theta$ is 
\begin{equation*}
\begin{split}
\Pi_\theta 
=\ &	\alpha_{01}tx \partial_t\wedge \partial_x +\alpha_{02}ty \partial_t\wedge \partial_y
	+\alpha_{03}tz \partial_t\wedge \partial_z \\ 
& 	+\alpha_{12}xy \partial_x\wedge \partial_y
	+\alpha_{13}xz \partial_x\wedge \partial_z +\alpha_{23}yz \partial_y\wedge \partial_z\,.
\end{split}
\end{equation*}
\item $A=\begin{pmatrix}a&1&&\\&a&&\\&&-a&1\\&&&-a\end{pmatrix}$, $a\neq0$. The compatible 1-form is 
\begin{equation*}\begin{split}
\theta =\	&  ty ( \alpha_1 y\,dt+\alpha_2t\,dy) +(tz-xy)(\delta_1(z\,dt-y\,dx)+\delta_2(t\,dz-x\,dy)) \\
	&  +ty( \beta_1(z\,dt-y\,dx)+\beta_2(t\,dz-x\,dy)) +(tz-xy)(\gamma_1y\,dt+\gamma_2t\,dy) 
\end{split}
\end{equation*}
and condition (ii) yields
\[
\delta_1-\delta_2=0,\qquad (\beta_1-\beta_2)( (\beta_1-\beta_2)+(\gamma_1+\gamma_2)=0\,.
\]
The first equation makes the second summand of $\theta$ a total derivative, such that it does not enter into the Poisson structure. A sample structure is given by the further choice $\beta_1=\beta_2$, $\gamma_1=\gamma_2$, i.e.\ 
\[
\Pi_\theta= \alpha ty \partial_x\wedge\partial_z  + \beta \Psi( d(ty)\wedge d(tz-xy) )\,.
\]
\item $A=\begin{pmatrix}a&d&&\\-d&a&&\\&&-a&e\\&&-e&-a\end{pmatrix}$, $de\neq0$, $d^2\neq e^2$

The compatible 1-form is 
\begin{equation*}
\begin{split}
\theta =\ & (y^2+z^2)(\alpha_1(x\,dt-t\,dx)+\alpha_2(t\,dt+x\,dx) ) \\
	& +(t^2+x^2)(\beta_1(z\,dy-y\,dz)+\beta_2(y\,dy+z\,dz) )
\end{split}
\end{equation*}
which also satisfies (ii). A sample structure is given by the choice $\alpha_2=\beta_2$ which combines the respective summands to a total derivative. The further choice $\alpha_1=\beta_1$ yields 
\begin{equation*}\begin{split}
\Pi_\theta=\ & (t^2+x^2)\partial_t\wedge\partial_x +(y^2+z^2)\partial_y\wedge\partial_z \\
	&+(ty+xz)(\partial_x\wedge\partial_y-\partial_t\wedge\partial_z)
	+(tz-xy)(\partial_t\wedge\partial_y+\partial_x\wedge\partial_z)
\end{split}\end{equation*}

\end{itemize}


\def\cprime{$'$} \def\cprime{$'$}

\end{document}